\documentclass[12pt]{article}

\usepackage{amsmath, epsfig, cite}
\usepackage{amssymb}
\usepackage{amsfonts}
\usepackage{latexsym}
\usepackage{amsthm}

\newtheorem{thm}{Theorem}[section]

\numberwithin{equation}{section}

\renewcommand{\thefootnote}

\begin{document}

\begin{center}
{\large\bf Double series for $\pi$ and their $q$-analogues
 \footnote{ Email addresses: weichuanan78@163.com (C. Wei), ythainmc@163.com (G. Ruan).}}
\end{center}

\renewcommand{\thefootnote}{$\dagger$}

\vskip 2mm \centerline{Chuanan Wei$^1$, Guozhu Ruan$^{2}$}
\begin{center}
{$^1$School of Biomedical Information and Engineering\\
  Hainan Medical University, Haikou 571199, China\\
  $^2$Medical Simulation Education Center\\
Hainan Medical University, Haikou 571199, China}
\end{center}


\vskip 0.7cm \noindent{\bf Abstract.} With the help of the partial
derivative operator and several summation formulas for
hypergeometric series, we find three double series for $\pi$. In
terms of the operator just stated and several summation formulas for
basic hypergeometric series, we also establish $q$-analogues of
these double series.



\section{Introduction}

For a complex variable $x$, define the well-known Gamma function to
be
\begin{align*}
\Gamma(x)=\int_{0}^{\infty}t^{x-1}e^{-t}dt\quad\text{with}\quad
Re(x)>0.
\end{align*}
Three important properties of this function can be expressed as
\begin{align*}
\Gamma(x+1)=x\Gamma(x),\quad
\Gamma(x)\Gamma(1-x)=\frac{\pi}{\sin(\pi x)}, \quad\lim_{n\to
\infty}\frac{\Gamma(x+n)}{\Gamma(y+n)}n^{y-x}=1,
\end{align*}
which will often be used without explanation in this paper.
Subsequently, we may give the definition of the shifted-factorial:
\begin{align*}
(x)_n=\Gamma(x+n)/\Gamma(x),
\end{align*}
where $n$ is an integer and $x$ is a complex number. Then the
hypergeometric series can be defined by
$$
_{r}F_{s}\left[\begin{array}{c}
a_1,a_2,\ldots,a_{r}\\
b_1,b_2,\ldots,b_{s}
\end{array};\, z
\right] =\sum_{k=0}^{\infty}\frac{(a_1)_k(a_2)_k\cdots(a_{r})_k}
{(b_1)_k(b_2)_k\cdots(b_{s})_k}\frac{z^k}{k!}.
$$

The research of $\pi$-formulas has a long history. In 1859, Bauer
\cite{Bauer} discovered a simple result:
\begin{align}
&\sum_{k=0}^{\infty}(-1)^k(4k+1)\frac{(\frac{1}{2})_k^3 }{k!^3 }
=\frac{2}{\pi}.
 \label{Bauer}
\end{align}
In 1914, Ramanujan \cite{Ramanujan} displayed 17 series for $1/\pi$
without proof. Decades later, Borweins \cite{Borwein} proved all of
them firstly. Two of Ramanujan's formulas are stated as
\begin{align}
&\sum_{k=0}^{\infty}(6k+1)\frac{(\frac{1}{2})_k^3 }{k!^3 4^k}
=\frac{4}{\pi},
 \label{Ramanujan-a}\\
&\sum_{k=0}^{\infty}(8k+1)\frac{(\frac{1}{2})_k(\frac{1}{4})_k(\frac{3}{4})_k}{k!^3
9^k} =\frac{2\sqrt{3}}{\pi}.
 \label{Ramanujan-b}
\end{align}

In 2011, Long \cite{Long} proposed the following conjecture : for
any odd prime $p$,
\begin{align}
\sum_{k=0}^{(p-1)/2}(-1)^k(6k+1)\frac{(\frac{1}{2})_k^3}{k!^3
8^k}\sum_{j=1}^{k}\left(\frac{1}{(2j-1)^2}-\frac{1}{16j^2}\right)
\equiv 0\pmod{p}, \label{eq:long-swisher}
\end{align}
which was certified by Swisher \cite{Swisher} after several years.
Recently, Guo and Lian \cite{Guo-c} conjectured two interesting
double series for $\pi$ related to \eqref{Ramanujan-a} and
\eqref{eq:long-swisher}:
\begin{align}
&\sum_{k=1}^{\infty}(6k+1)\frac{(\frac{1}{2})_k^3 }{k!^3
4^k}\sum_{j=1}^{k}\bigg\{\frac{1}{(2j-1)^2}-\frac{1}{16j^2}\bigg\}
=\frac{\pi}{12}, \label{eq:guo-a}
\\
 &\sum_{k=1}^{\infty}(-1)^k(6k+1)\frac{(\frac{1}{2})_k^3}{k!^3
8^k}\sum_{j=1}^{k}\bigg\{\frac{1}{(2j-1)^2}-\frac{1}{16j^2}\bigg\}
=-\frac{\sqrt{2}}{48}\pi.
 \label{eq:guo-b}
\end{align}
which have been proved by Wei \cite{Wei}.
 For more
known series on $\pi$, we refer the reader to the papers
\cite{Guillera-b,Schlosser,Chan,Liu,Sun-a,Wang,Zudilin}.

Inspired by the work just mentioned, we shall established the
following two theorems associated with \eqref{Bauer} and
\eqref{Ramanujan-b}.

\begin{thm}\label{thm-a}
\begin{align}
\sum_{k=1}^{\infty}(-1)^k(4k+1)\frac{(\frac{1}{2})_k^3 }{k!^3}
\sum_{i=1}^{2k}\frac{(-1)^i}{i^2} =\frac{\pi}{12}.
 \label{eq:wei-a}
\end{align}
\end{thm}

\begin{thm}\label{thm-b}
\begin{align}
&\sum_{k=1}^{\infty}(8k+1)\frac{(\frac{1}{2})_k(\frac{1}{4})_k(\frac{3}{4})_k}{k!^3\,
9^k} \sum_{i=1}^{k}\bigg\{\frac{1}{(2i-1)^2}-\frac{1}{36i^2}\bigg\}
=\frac{\sqrt{3}\,\pi}{54}.
 \label{eq:wei-b}
\end{align}
\end{thm}

Furthermore, we shall provide the following double series for
$\pi^3$.

\begin{thm}\label{thm-c}
\begin{align}
&\sum_{k=1}^{\infty}(8k-1)\frac{(1)_k(-\frac{1}{4})_k(-\frac{3}{4})_k}{(\frac{1}{2})_k^2(\frac{3}{2})_k\,9^k}
\sum_{i=1}^{k}\bigg\{\frac{1}{(2i-1)^2}-\frac{9}{4i^2}\bigg\}
=-\frac{\sqrt{3}\,\pi^3}{108}.
 \label{eq:wei-c}
\end{align}
\end{thm}

For an integer $n$ and two complex numbers $x$, $q$ with $|q|<1$,
define the $q$-shifted factorial to be
\begin{align*}
(x;q)_{\infty}=\prod_{i=0}^{\infty}(1-xq^i),\quad
(x;q)_n=\frac{(x;q)_{\infty}}{(xq^n;q)_{\infty}}.
\end{align*}
For convenience, we shall also adopt  the following notation:
\begin{equation*}
(x_1,x_2,\dots,x_r;q)_{m}=(x_1;q)_{m}(x_2;q)_{m}\cdots(x_r;q)_{m},
 \end{equation*}
where $r\in\mathbb{Z}^{+}$ and $m\in\mathbb{Z}^{+}\cup\{0,\infty\}$.
Then following Gasper and Rahman \cite{Gasper}, the basic
hypergeometric series can be defined as
$$
_{r}\phi_{s}\left[\begin{array}{c}
a_1,a_2,\ldots,a_{r}\\
b_1,b_2,\ldots,b_{s}
\end{array};q,\, z
\right] =\sum_{k=0}^{\infty}\frac{(a_1,a_2,\ldots, a_{r};q)_k}
{(q,b_1,b_2,\ldots,b_{s};q)_k}\bigg\{(-1)^kq^{\binom{k}{2}}\bigg\}^{1+s-r}z^k.
$$

Let $[n]=1+q+\cdots+q^{n-1}$ be the $q$-integer. Recently, Guo and
Liu \cite{Guo-d} and Guo and Zudilin \cite{Guo-e} obtained the
$q$-analogues of \eqref{Bauer}-\eqref{Ramanujan-b}:
 \begin{align*}
&\sum_{k=0}^{\infty}(-1)^kq^{k^2}[4k+1]\frac{(q;q^2)_k^3}{(q^2;q^2)_k^3}
=\frac{(q,q^3;q^2)_{\infty}}{(q^2;q^2)_{\infty}^2},
 \\
 &\sum_{k=0}^{\infty}q^{k^2}[6k+1]\frac{(q;q^2)_k^2(q^2;q^4)_k}{(q^4;q^4)_k^3}
=\frac{(1+q)(q^2,q^6;q^4)_{\infty}}{(q^4;q^4)_{\infty}^2},
\\
&\sum_{k=0}^{\infty}q^{2k^2}[8k+1]\frac{(q;q^2)_{k}^2(q;q^2)_{2k}}{(q^2;q^2)_{2k}(q^6;q^6)_{k}^2}
=\frac{(q^3;q^2)_{\infty}(q^3;q^6)_{\infty}}{(q^2;q^2)_{\infty}(q^6;q^6)_{\infty}}.
 \end{align*}
Wei \cite{Wei} got the $q$-analogues of \eqref{eq:guo-a} and
\eqref{eq:guo-b}:
\begin{align*}
&\sum_{k=1}^{\infty}q^{k^2}[6k+1]\frac{(q;q^2)_k^2(q^2;q^4)_k}{(q^4;q^4)_k^3}
\sum_{j=1}^{k}\bigg\{\frac{q^{2j-1}}{[2j-1]^2}-\frac{q^{4j}}{[4j]^2}\bigg\}
\\[3pt]
&\quad=\frac{(q^2;q^4)_{\infty}^2(q^5;q^4)_{\infty}}{(q;q^4)_{\infty}(q^4;q^4)_{\infty}^2}
\sum_{i=1}^{\infty}(-1)^{i-1}\frac{q^{2i}}{[2i]^2},
\\[3pt]
&\sum_{k=1}^{\infty}(-1)^kq^{3k^2}[6k+1]\frac{(q;q^2)_k^3}{(q^4;q^4)_k^3}
\sum_{j=1}^{k}\bigg\{\frac{q^{2j-1}}{[2j-1]^2}-\frac{q^{4j}}{[4j]^2}\bigg\}
\\[3pt]
&\quad=-\frac{(q^3,q^5;q^4)_{\infty}}{(q^4;q^4)_{\infty}^2}
\sum_{i=1}^{\infty}\frac{q^{4i}}{[4i]^2}.
\end{align*}
More $q$-analogues of
 $\pi$-formulas can be seen in the papers
 \cite{Guo-b,Guo-f,He,Hou-a,Hou-b,Sun-b}.

 Inspired by the work
 just mentioned, we shall derive the following $q$-analogues of
 Theorems \ref{thm-a}-\ref{thm-c}.

\begin{thm}\label{thm-d}
\begin{align}
\sum_{k=1}^{\infty}(-1)^kq^{k^2}[4k+1]\frac{(q;q^2)_k^3}{(q^2;q^2)_k^3}\sum_{i=1}^{2k}(-1)^i\frac{q^i}{[i]^2}
=\frac{(q,q^3;q^2)_{\infty}}{(q^2;q^2)_{\infty}^2}\sum_{j=1}^{\infty}\frac{q^{2j}}{[2j]^2}.
\label{eq:wei-d}
\end{align}
\end{thm}

\begin{thm}\label{thm-e}
\begin{align}
&\sum_{k=1}^{\infty}q^{2k^2}[8k+1]\frac{(q;q^2)_{k}^2(q;q^2)_{2k}}{(q^2;q^2)_{2k}(q^6;q^6)_{k}^2}
\sum_{i=1}^{k}\bigg\{\frac{q^{2i-1}}{[2i-1]^2}-\frac{q^{6i}}{[6i]^2}\bigg\}
\notag\\[3pt]
&\quad=\frac{(q^3;q^2)_{\infty}(q^3;q^6)_{\infty}}{(q^2;q^2)_{\infty}(q^6;q^6)_{\infty}}
\sum_{j=1}^{\infty}(-1)^{j-1}\frac{q^{3j}}{[3j]^2}.
  \label{eq:wei-e}
\end{align}
\end{thm}

\begin{thm}\label{thm-f}
\begin{align}
&\sum_{k=1}^{\infty}q^{2k^2+2k}[8k-1]\frac{(q^2;q^2)_{k}^2(q^{-3};q^2)_{2k}}{(q^4;q^2)_{2k}(q^3;q^6)_{k}^2}
\sum_{i=1}^{k}\bigg\{\frac{q^{6i-3}}{[6i-3]^2}-\frac{q^{2i}}{[2i]^2}\bigg\}
\notag\\[3pt]
&\quad=\frac{(q;q^2)_{\infty}(q^6;q^6)_{\infty}^3}{(q^4;q^2)_{\infty}(q^3;q^6)_{\infty}^3}
\sum_{j=1}^{\infty}(-1)^{j}\frac{q^{3j-1}}{[3j]^2}.
 \label{eq:wei-f}
\end{align}
\end{thm}

For a multivariable function $f(x_1,x_2,\ldots,x_m)$, define the
partial derivative operator $\mathcal{D}_{x_i}$ by
\begin{align*}
&&\mathcal{D}_{x_i}f(x_1,x_2,\ldots,x_m)=\frac{d}{dx_i}f(x_1,x_2,\ldots,x_m)\quad\text{with}\quad
1\leq i\leq m.
 \end{align*}
Then we have the following four relations: for $n>0$,
\begin{align*}
&\mathcal{D}_{x}(x+y)_n=(x+y)_n\sum_{i=1}^n\frac{1}{x+y-1+i},
\\[3pt]
&\mathcal{D}_{x}\sum_{i=1}^n\frac{1}{x+y+i}=-\sum_{i=1}^n\frac{1}{(x+y+i)^2},
\\[3pt]
&\mathcal{D}_{x}(xy;q)_n=-(xy;q)_n\sum_{i=1}^n\frac{yq^{i-1}}{1-xyq^{i-1}},
\\[3pt]
&\mathcal{D}_{x}\sum_{i=1}^n\frac{yq^{i}}{1-xyq^{i}}=\sum_{i=1}^n\frac{y^2q^{2i}}{(1-xyq^{i})^2},
\end{align*}
which will frequently be utilized without indication in this paper.

The structure of the paper is arranged as follows. We shall verify
Theorems \ref{thm-a}-\ref{thm-c} via the partial derivative operator
and some summation formulas for hypergeometric series in Section 2.
Similarly, we shall prove Theorems \ref{thm-d}-\ref{thm-f} through
the partial derivative operator and some summation formulas for
basic hypergeometric series in Section 3.

\section{Proof of Theorems \ref{thm-a}-\ref{thm-c}}

Firstly, we shall prove Theorem \ref{thm-a}.

\begin{proof}[{\bf{Proof of Theorem \ref{thm-a}}}]
 Recall Dougall's $_5F_4$ summation formula (cf. \cite[P.
71]{Andrews}):
\begin{align*}
{_{5}F_{4}}\left[\begin{array}{cccccccc}
  a,1+\frac{a}{2},b,c,-n\\
  \frac{a}{2},1+a-b,1+a-c,1+a+n
\end{array};1
\right]=\frac{(1+a)_n(1+a-b-c)_n}{(1+a-b)_n(1+a-c)_n}.
\end{align*}
The $c=1-b$ case of it reads
\begin{align}
{_{5}F_{4}}\left[\begin{array}{cccccccc}
  a,1+\frac{a}{2},b,1-b,-n\\
  \frac{a}{2},1+a-b,a+b,1+a+n
\end{array};1
\right]=\frac{(a)_n(1+a)_n}{(a+b)_n(1+a-b)_n}.
 \label{eq:wei-aa}
\end{align}
Apply the partial derivative operator $\mathcal{D}_{b}$ to both
sides of \eqref{eq:wei-aa} to obtain
\begin{align*}
&\sum_{k=1}^n \frac{(a)_k(1+\frac{a}{2})_k(b)_k(1-b)_k(-n)_k}
{(1)_k(\frac{a}{2})_k(1+a-b)_k(a+b)_k(1+a+n)_k}
\notag\\[2pt]
&\quad\times\Bigg\{\sum_{i=1}^k\frac{1}{b-1+i}-\sum_{i=1}^k\frac{1}{-b+i}+\sum_{i=1}^k\frac{1}{a-b+i}-\sum_{i=1}^k\frac{1}{a+b-1+i}\Bigg\}
\notag\\[2pt]
&\:=\frac{(a)_n(1+a)_n}{(a+b)_n(1+a-b)_n}
\Bigg\{\sum_{j=1}^n\frac{1}{a-b+j}-\sum_{j=1}^n\frac{1}{a+b-1+j}\Bigg\}.
\end{align*}
Employing the operator $\mathcal{D}_{b}$ to both sides of the last
equation, there holds
\begin{align}
&\sum_{k=1}^n \frac{(a)_k(1+\frac{a}{2})_k(b)_k(1-b)_k(-n)_k}
{(1)_k(\frac{a}{2})_k(1+a-b)_k(a+b)_k(1+a+n)_k}
\notag\\[2pt]
&\quad\,\,\times\Bigg\{\bigg[\sum_{i=1}^k\frac{1}{b-1+i}-\sum_{i=1}^k\frac{1}{-b+i}+\sum_{i=1}^k\frac{1}{a-b+i}-\sum_{i=1}^k\frac{1}{a+b-1+i}\bigg]^2
\notag\\[2pt]
&\qquad\:-\bigg[\sum_{i=1}^k\frac{1}{(b-1+i)^2}+\sum_{i=1}^k\frac{1}{(-b+i)^2}-\sum_{i=1}^k\frac{1}{(a-b+i)^2}-\sum_{i=1}^k\frac{1}{(a+b-1+i)^2}\bigg]\Bigg\}
\notag\\[2pt]
&\:=\frac{(a)_n(1+a)_n}{(a+b)_n(1+a-b)_n}
\Bigg\{\bigg[\sum_{j=1}^n\frac{1}{a-b+j}-\sum_{j=1}^n\frac{1}{a+b-1+j}\bigg]^2
\notag\\[2pt]
&\qquad\:+\bigg[\sum_{j=1}^n\frac{1}{(a-b+j)^2}+\sum_{j=1}^n\frac{1}{(a+b-1+j)^2}\bigg]\Bigg\}.
\label{eq:wei-bb}
\end{align}
The $a=b=\frac{1}{2}$ case of \eqref{eq:wei-bb} can be manipulated
as
\begin{align*}
&\sum_{k=1}^{n}(4k+1)\frac{(\frac{1}{2})_k^3 }{k!^3}
\frac{(-n)_k}{(\frac{3}{2}+n)_k}\sum_{i=1}^{2k}\frac{(-1)^i}{i^2}
=\frac{\Gamma(\frac{1}{2}+n)\Gamma(\frac{3}{2}+n)}{\Gamma(1+n)\Gamma(1+n)}
\frac{1}{\Gamma(\frac{1}{2})\Gamma(\frac{3}{2})}\sum_{j=1}^{n}\frac{1}{4j^2}.
\end{align*}
Letting $n\to \infty$ in the above identity,
 we arrive at
\begin{align*}
\sum_{k=1}^{\infty}(-1)^k(4k+1)\frac{(\frac{1}{2})_k^3 }{k!^3}
\sum_{i=1}^{2k}\frac{(-1)^i}{i^2}
=\frac{1}{\Gamma(\frac{1}{2})\Gamma(\frac{3}{2})}\sum_{j=1}^{\infty}\frac{1}{4j^2}.
\end{align*}
Calculating the series, which is on the right-hand side, by  Euler's
formula:
\begin{align}\label{Euler}
\sum_{j=1}^{\infty}\frac{1}{j^2}=\frac{\pi^2}{6},
\end{align}
we are led to \eqref{eq:wei-a}.

\end{proof}

Secondly, we shall give the proof of Theorem \ref{thm-b}.

\begin{proof}[{\bf{Proof of Theorem \ref{thm-b}}}]
 An original Gosper Conjecture (cf. \cite[p. 307]{GD}) is
\begin{align}
&{_7}F_{6}\left[\begin{array}{c}
a,1+\frac{a}{2},a+\frac{1}{2},b,1-b,\frac{2a+1}{3}+n,-n\\[3pt]
\frac{a}{2},\frac{1}{2},\frac{2a-b+3}{3},\frac{2a+b+2}{3},-3n,1+2a+3n
\end{array};\, 1 \right]
=\frac{(\frac{1+b}{3})_n(\frac{2-b}{3})_n(\frac{2a+2}{3})_n(\frac{2a+3}{3})_n}
{(\frac{1}{3})_n(\frac{2}{3})_n(\frac{2a-b+3}{3})_n(\frac{2a+b+2}{3})_n},
\label{Gosper}
\end{align}
the nonterminating form of which can be seen in Gasper and Rahman
\cite[Equation (1.6)]{Gasper-Rahman}. By means of the partial
derivative operator $\mathcal{D}_{b}$ and the $a\to a/2$ case of
\eqref{Gosper}, we have
\begin{align*}
&\sum_{k=1}^n
\frac{(\frac{a}{2})_k(1+\frac{a}{4})_k(\frac{a+1}{2})_k(b)_k(1-b)_k(\frac{a+1}{3}+n)_k(-n)_k}
{(1)_k(\frac{a}{4})_k(\frac{1}{2})_k(\frac{a-b+3}{3})_k(\frac{a+b+2}{3})_k(-3n)_k(1+a+3n)_k}
\notag\\[2pt]
&\quad\times\Bigg\{\sum_{i=1}^k\frac{3}{b-1+i}-\sum_{i=1}^k\frac{3}{-b+i}+\sum_{i=1}^k\frac{1}{\frac{a-b}{3}+i}-\sum_{i=1}^k\frac{1}{\frac{a+b-1}{3}+i}\Bigg\}
\notag\\[2pt]
&=\frac{(\frac{1+b}{3})_n(\frac{2-b}{3})_n(\frac{a+2}{3})_n(\frac{a+3}{3})_n}
{(\frac{1}{3})_n(\frac{2}{3})_n(\frac{a-b+3}{3})_n(\frac{a+b+2}{3})_n}
\notag\\[2pt]
&\quad\times\Bigg\{\sum_{j=1}^n\frac{1}{\frac{b-2}{3}+j}-\sum_{j=1}^n\frac{1}{\frac{-b-1}{3}+j}+\sum_{j=1}^n\frac{1}{\frac{a-b}{3}+j}-\sum_{j=1}^n\frac{1}{\frac{a+b-1}{3}+j}\Bigg\}.
\end{align*}
According to the operator $\mathcal{D}_{b}$ and the last equation,
it is routine to understand that
\begin{align}
&\sum_{k=1}^n
\frac{(\frac{a}{2})_k(1+\frac{a}{4})_k(\frac{a+1}{2})_k(b)_k(1-b)_k(\frac{a+1}{3}+n)_k(-n)_k}
{(1)_k(\frac{a}{4})_k(\frac{1}{2})_k(\frac{a-b+3}{3})_k(\frac{a+b+2}{3})_k(-3n)_k(1+a+3n)_k}
\notag\\[2pt]
&\quad\times\Bigg\{\bigg[\sum_{i=1}^k\frac{3}{b-1+i}-\sum_{i=1}^k\frac{3}{-b+i}
+\sum_{i=1}^k\frac{1}{\frac{a-b}{3}+i}-\sum_{i=1}^k\frac{1}{\frac{a+b-1}{3}+i}\bigg]^2
\notag\\[2pt]
&\qquad-\bigg[\sum_{i=1}^k\frac{9}{(b-1+i)^2}+\sum_{i=1}^k\frac{9}{(-b+i)^2}
-\sum_{i=1}^k\frac{1}{(\frac{a-b}{3}+i)^2}-\sum_{i=1}^k\frac{1}{(\frac{a+b-1}{3}+i)^2}\bigg]\Bigg\}
\notag
\end{align}
\begin{align}
&=\frac{(\frac{1+b}{3})_n(\frac{2-b}{3})_n(\frac{a+2}{3})_n(\frac{a+3}{3})_n}
{(\frac{1}{3})_n(\frac{2}{3})_n(\frac{a-b+3}{3})_n(\frac{a+b+2}{3})_n}
\notag\\[2pt]
&\quad\times\Bigg\{\bigg[\sum_{j=1}^n\frac{1}{\frac{b-2}{3}+j}-\sum_{i=1}^n\frac{1}{\frac{-b-1}{3}+j}
+\sum_{j=1}^n\frac{1}{\frac{a-b}{3}+j}-\sum_{j=1}^n\frac{1}{\frac{a+b-1}{3}+j}\bigg]^2
\notag\\[2pt]
&\qquad-\bigg[\sum_{j=1}^n\frac{1}{(\frac{b-2}{3}+j)^2}+\sum_{j=1}^n\frac{1}{(\frac{-b-1}{3}+j)^2}
-\sum_{j=1}^n\frac{1}{(\frac{a-b}{3}+j)^2}-\sum_{j=1}^n\frac{1}{(\frac{a+b-1}{3}+j)^2}\bigg]\Bigg\}.
\label{eq:wei-cc}
\end{align}
The $a=b=\frac{1}{2}$ case of \eqref{eq:wei-cc} produces
\begin{align*}
&\sum_{k=1}^{n}(8k+1)\frac{(\frac{1}{2})_k(\frac{1}{4})_k(\frac{3}{4})_k}{k!^3}
\frac{(\frac{1}{2}+n)_k(-n)_k}{(\frac{3}{2}+3n)_k(-3n)_k}\sum_{i=1}^{k}\bigg\{\frac{1}{(2i-1)^2}-\frac{1}{36i^2}\bigg\}
\\[2pt]
&\quad=\frac{1}{9}\frac{\Gamma(\frac{1}{2}+n)^2\Gamma(\frac{5}{6}+n)\Gamma(\frac{7}{6}+n)}
{\Gamma(1+n)^2\Gamma(\frac{1}{3}+n)\Gamma(\frac{2}{3}+n)}
\frac{\Gamma(\frac{1}{3})\Gamma(\frac{2}{3})}{\Gamma(\frac{1}{2})^2\Gamma(\frac{5}{6})\Gamma(\frac{7}{6})}\sum_{j=1}^{2n}\frac{(-1)^{j-1}}{j^2}.
\end{align*}
Letting $n\to \infty$ in the upper identity, we get
\begin{align}
&\sum_{k=1}^{\infty}(8k+1))\frac{(\frac{1}{2})_k(\frac{1}{4})_k(\frac{3}{4})_k}{k!^3\,9^k}
\sum_{j=1}^{k}\bigg\{\frac{1}{(2j-1)^2}-\frac{1}{16j^2}\bigg\}
\notag\\[2mm]
&\quad=\frac{1}{9}\frac{\Gamma(\frac{1}{3})\Gamma(\frac{2}{3})}{\Gamma(\frac{1}{2})^2\Gamma(\frac{5}{6})\Gamma(\frac{7}{6})}\sum_{j=1}^{\infty}\frac{(-1)^{j-1}}{j^2}.
\label{eq:wei-dd}
\end{align}
 On the base of Euler's formula \eqref{Euler},
we can find  that
\begin{align*}
\sum_{j=1}^{\infty}\frac{1}{(2j-1)^2}=\sum_{j=1}^{\infty}\frac{1}{j^2}-\sum_{j=1}^{\infty}\frac{1}{(2j)^2}=\frac{\pi^2}{8}.
\end{align*}
So there is the formula
\begin{align}
\sum_{j=1}^{\infty}\frac{(-1)^{j-1}}{j^2}=\sum_{j=1}^{\infty}\frac{1}{(2j-1)^2}-\sum_{i=1}^{\infty}\frac{1}{(2j)^2}=\frac{\pi^2}{12}.
\label{eq:wei-ee}
\end{align}
Substituting \eqref{eq:wei-ee} into \eqref{eq:wei-dd}, we catch hold
of \eqref{eq:wei-b}.
\end{proof}

Thirdly, we shall display the proof of Theorem \ref{thm-c}.

\begin{proof}[{\bf{Proof of Theorem \ref{thm-c}}}]
 A known $_7F_6$ summation formula ((cf. \cite[p. 37]{Chu}) can be
 written as
\begin{align}
&{_7}F_{6}\left[\begin{array}{c}
a,1+\frac{a}{2},a-\frac{1}{2},b,2-b,\frac{2a+2}{3}+n,-n\\[3pt]
\frac{a}{2},\frac{3}{2},\frac{2a-b+3}{3},\frac{2a+b+1}{3},-1-3n,1+2a+3n
\end{array};\, 1 \right]
\notag\\[3pt]
  &\quad
=\frac{(\frac{2+b}{3})_n(\frac{4-b}{3})_n(\frac{2a+1}{3})_n(\frac{2a+3}{3})_n}
{(\frac{2}{3})_n(\frac{4}{3})_n(\frac{2a-b+3}{3})_n(\frac{2a+b+1}{3})_n},
\label{equation-a}
\end{align}
the nonterminating form of which can be observed in Gasper and
Rahman \cite[Equation (4.7)]{Gasper-Rahman}. Apply the partial
derivative operator $\mathcal{D}_{b}$ to the $a\to a/2$ case of
\eqref{equation-a} to deduce
\begin{align*}
&\sum_{k=1}^n
\frac{(\frac{a}{2})_k(1+\frac{a}{4})_k(\frac{a-1}{2})_k(b)_k(2-b)_k(\frac{a+2}{3}+n)_k(-n)_k}
{(1)_k(\frac{a}{4})_k(\frac{3}{2})_k(\frac{a-b+3}{3})_k(\frac{a+b+1}{3})_k(-1-3n)_k(1+a+3n)_k}
\notag\\[2pt]
&\quad\times\Bigg\{\sum_{i=1}^k\frac{3}{b-1+i}-\sum_{i=1}^k\frac{3}{1-b+i}+\sum_{i=1}^k\frac{1}{\frac{a-b}{3}+i}-\sum_{i=1}^k\frac{1}{\frac{a+b-2}{3}+i}\Bigg\}
\notag\\[2pt]
&=\frac{(\frac{2+b}{3})_n(\frac{4-b}{3})_n(\frac{a+1}{3})_n(\frac{a+3}{3})_n}
{(\frac{2}{3})_n(\frac{4}{3})_n(\frac{a-b+3}{3})_n(\frac{a+b+1}{3})_n}
\notag\\[2pt]
&\quad\times\Bigg\{\sum_{j=1}^n\frac{1}{\frac{b-1}{3}+j}-\sum_{j=1}^n\frac{1}{\frac{1-b}{3}+j}+\sum_{j=1}^n\frac{1}{\frac{a-b}{3}+j}-\sum_{j=1}^n\frac{1}{\frac{a+b-2}{3}+j}\Bigg\}.
\end{align*}
Employing the operator $\mathcal{D}_{b}$ to both sides of the last
equation, it is easy to realize that
\begin{align}
&\sum_{k=1}^n
\frac{(\frac{a}{2})_k(1+\frac{a}{4})_k(\frac{a-1}{2})_k(b)_k(2-b)_k(\frac{a+2}{3}+n)_k(-n)_k}
{(1)_k(\frac{a}{4})_k(\frac{3}{2})_k(\frac{a-b+3}{3})_k(\frac{a+b+1}{3})_k(-1-3n)_k(1+a+3n)_k}
\notag\\[2pt]
&\quad\times\Bigg\{\bigg[\sum_{i=1}^k\frac{3}{b-1+i}-\sum_{i=1}^k\frac{3}{1-b+i}+\sum_{i=1}^k\frac{1}{\frac{a-b}{3}+i}-\sum_{i=1}^k\frac{1}{\frac{a+b-2}{3}+i}\bigg]^2
\notag\\[2pt]
&\qquad-\bigg[\sum_{i=1}^k\frac{9}{(b-1+i)^2}+\sum_{i=1}^k\frac{9}{(1-b+i)^2}
-\sum_{i=1}^k\frac{1}{(\frac{a-b}{3}+i)^2}-\sum_{i=1}^k\frac{1}{(\frac{a+b-2}{3}+i)^2}\bigg]\Bigg\}
\notag\\[2pt]
&=\frac{(\frac{2+b}{3})_n(\frac{4-b}{3})_n(\frac{a+1}{3})_n(\frac{a+3}{3})_n}
{(\frac{2}{3})_n(\frac{4}{3})_n(\frac{a-b+3}{3})_n(\frac{a+b+1}{3})_n}
\notag\\[2pt]
&\quad\times\Bigg\{\bigg[\sum_{j=1}^n\frac{1}{\frac{b-1}{3}+j}-\sum_{j=1}^n\frac{1}{\frac{1-b}{3}+j}
+\sum_{j=1}^n\frac{1}{\frac{a-b}{3}+j}-\sum_{j=1}^n\frac{1}{\frac{a+b-2}{3}+j}\bigg]^2
\notag\\[2pt]
&\qquad-\bigg[\sum_{j=1}^n\frac{1}{(\frac{b-1}{3}+j)^2}+\sum_{j=1}^n\frac{1}{(\frac{1-b}{3}+j)^2}
-\sum_{j=1}^n\frac{1}{(\frac{a-b}{3}+j)^2}-\sum_{j=1}^n\frac{1}{(\frac{a+b-2}{3}+j)^2}\bigg]\Bigg\}.
\label{eq:wei-hh}
\end{align}
The $a=-\frac{1}{2}, b=1$ case of \eqref{eq:wei-hh} engenders
\begin{align*}
&\sum_{k=1}^{n}(8k-1)\frac{(1)_k(-\frac{1}{4})_k(-\frac{3}{4})_k}{(\frac{1}{2})_k^2(\frac{3}{2})_k}
\frac{(\frac{1}{2}+n)_k(-n)_k}{(\frac{1}{2}+3n)_k(-1-3n)_k}\sum_{i=1}^{k}\bigg\{\frac{1}{(2i-1)^2}-\frac{9}{4i^2}\bigg\}
\\[2pt]
&\quad=\frac{\Gamma(1+n)^2\Gamma(\frac{1}{6}+n)\Gamma(\frac{5}{6}+n)}
{\Gamma(\frac{1}{2}+n)^2\Gamma(\frac{2}{3}+n)\Gamma(\frac{4}{3}+n)}
\frac{\Gamma(\frac{1}{2})^2\Gamma(\frac{2}{3})\Gamma(\frac{4}{3})}{\Gamma(\frac{1}{6})\Gamma(\frac{5}{6})}\sum_{j=1}^{2n}\frac{(-1)^{j}}{j^2}.
\end{align*}
Letting $n\to \infty$ in this identity and using \eqref{eq:wei-ee},
we discover \eqref{eq:wei-c}.
\end{proof}

\section{Proof of Theorems \ref{thm-d}-\ref{thm-f}}

For proving Theorem \ref{thm-d}, we need the $q$-analogue of
Dougall's $_5F_4$ summation formula (cf. \cite[Equation
(2.4.2)]{Gasper}):
\begin{align}
{_6}\phi_{5}\left[\begin{array}{c}
a,qa^{\frac{1}{2}},-qa^{\frac{1}{2}},b,c,q^{-n}\\[3pt]
a^{\frac{1}{2}},-a^{\frac{1}{2}},aq/b,aq/c,aq^{n+1}
\end{array};\, q, \frac{aq^{n+1}}{bc} \right]
=\frac{(aq,aq/bc;q)_{n}}{(aq/b,aq/c;q)_{n}}. \label{basic-a}
\end{align}

Now we are ready to prove Theorem \ref{thm-d}.

\begin{proof}[{\bf{Proof of Theorem \ref{thm-d}}}]
Apply the partial derivative operator $\mathcal{D}_{b}$ to the
$c=q/b$ case of \eqref{basic-a} to obtain
\begin{align*}
&\sum_{k=1}^{n}\frac{1-aq^{2k}}{1-a}\frac{(a,b,q/b,q^{-n};q)_k}{(q,aq/b,ab,aq^{n+1};q)_k}(aq^n)^kA_k(a,b)
= \frac{(a,aq;q)_{n}}{(ab,aq/b;q)_{n}}B_n(a,b),
\end{align*}
where
\begin{align*}
&A_k(a,b)=\sum_{i=1}^k\frac{q^{i-1}}{1-bq^{i-1}}-\sum_{i=1}^k\frac{q^{i}/b^2}{1-q^{i}/b}
+\sum_{i=1}^k\frac{aq^{i}/b^2}{1-aq^{i}/b}-\sum_{i=1}^k\frac{aq^{i-1}}{1-abq^{i-1}},
\\[2pt]
&B_n(a,b)=\sum_{j=1}^n\frac{aq^{j}/b^2}{1-aq^{j}/b}-\sum_{j=1}^n\frac{aq^{j-1}}{1-abq^{j-1}}.
\end{align*}
Employing the operator $\mathcal{D}_{b}$ to both sides of the last
equation, there holds
\begin{align}
&\sum_{k=1}^{n}\frac{1-aq^{2k}}{1-a}\frac{(a,b,q/b,q^{-n};q)_k}{(q,aq/b,ab,aq^{n+1};q)_k}(aq^n)^k
\bigg\{A_k(a,b)^2-C_k(a,b)\bigg\}
\notag\\[3pt]
  &\quad
=\frac{(a,aq;q)_{n}}{(ab,aq/b;q)_{n}}
\bigg\{B_n(a,b)^2-D_n(a,b)\bigg\}, \label{eq:wei-aaa}
\end{align}
where
\begin{align*}
&C_k(a,b)=\sum_{i=1}^k\frac{q^{2i-2}}{(1-bq^{i-1})^2}-\sum_{i=1}^k\frac{(q^{i}/b-2)q^i/b^3}{(1-q^{i}/b)^2}
\\[2pt]
&\qquad\qquad\:+\sum_{i=1}^k\frac{(aq^{i}/b-2)aq^{i}/b^3}{(1-aq^{i}/b)^2}-\sum_{i=1}^k\frac{a^2q^{2i-2}}{(1-abq^{i-1})^2},
\\[2pt]
&D_n(a,b)=\sum_{j=1}^n\frac{(aq^{j}/b-2)aq^{j}/b^3}{(1-aq^{j}/b)^2}-\sum_{j=1}^n\frac{a^2q^{2j-2}}{(1-abq^{j-1})^2}.
\end{align*}
The $a\to q, b\to q, q\to q^2$ case of \eqref{eq:wei-aaa} reads
\begin{align*}
&\sum_{k=1}^{n}[4k+1]
\frac{(q;q^2)_k^3}{(q^2;q^2)_k^3}\frac{(q^{-2n};q^2)_k}{(q^{3+2n};q^2)_k}\,q^{(1+2n)k}
\sum_{i=1}^{2k}(-1)^i\frac{q^i}{[i]^2}
=\frac{(q,q^3;q^2)_{n}}{(q^2;q^2)_{n}^2}
\sum_{j=1}^n\frac{q^{2j}}{[2j]^2}.
\end{align*}
Letting $n\to \infty$ in the above identity, we arrive at
\eqref{eq:wei-d}.
\end{proof}

In order to prove Theorem \ref{thm-e}, we require the summation
formula for basic hypergeometric series  (cf. \cite[p. 65]{Chu-a}):
\begin{align}
&\sum_{k=0}^{n}\frac{1-a^{\frac{2}{3}}q^{\frac{4k}{3}}}{1-a^{\frac{2}{3}}}\frac{(a^{\frac{2}{3}};q^{\frac{1}{3}})_{2k}}{(q^{\frac{1}{3}};q^{\frac{1}{3}})_{2k}}
\frac{(b^{\frac{1}{3}},q^{\frac{1}{3}}b^{-\frac{1}{3}};q^{\frac{1}{3}})_{k}}{(a^{\frac{2}{3}}q^{n+\frac{1}{3}},q^{-n};q^{\frac{1}{3}})_{k}}
\frac{(a^{\frac{2}{3}}q^{n+\frac{1}{3}},q^{-n};q)_k}{(qa^{\frac{2}{3}}b^{-\frac{1}{3}},q^{\frac{2}{3}}a^{\frac{2}{3}}b^{\frac{1}{3}};q)_k}
q^{\frac{k}{3}}
\notag\\[3pt]
  &\quad
=\frac{(q^{\frac{1}{3}}a^{\frac{2}{3}};q^{\frac{1}{3}})_{n}}{(q^{\frac{1}{3}};q^{\frac{1}{3}})_{n}}\frac{(q,q^{\frac{1}{3}}b^{\frac{1}{3}},q^{\frac{2}{3}}b^{-\frac{1}{3}};q)_n}
{(q^{\frac{1}{3}}a^{\frac{2}{3}},qa^{\frac{2}{3}}b^{-\frac{1}{3}},q^{\frac{2}{3}}a^{\frac{2}{3}}b^{\frac{1}{3}};q)_n},
\label{basic-b}
\end{align}
where we have replaced
$$\frac{(b^{\frac{1}{3}},q^{\frac{1}{3}}b^{-\frac{1}{3}};q)_{k}}{(a^{\frac{2}{3}}q^{n+\frac{1}{3}},q^{-n};q)_{k}}
\quad\text{by}\quad
\frac{(b^{\frac{1}{3}},q^{\frac{1}{3}}b^{-\frac{1}{3}};q^{\frac{1}{3}})_{k}}{(a^{\frac{2}{3}}q^{n+\frac{1}{3}},q^{-n};q^{\frac{1}{3}})_{k}}$$
for correction. The nonterminating form of \eqref{basic-b} can be
seen in Gasper and Rahman \cite[Equation (1.8)]{Gasper-Rahman}.

Subsequently, we start to prove Theorem \ref{thm-e}.

\begin{proof}[{\bf{Proof of Theorem \ref{thm-e}}}]
Via the partial derivative operator $\mathcal{D}_{b}$ and the $a\to
a^{3/2}, b\to b^3, q\to q^3$ of \eqref{basic-b}, we get
\begin{align*}
&\sum_{k=1}^{n}\frac{1-aq^{4k}}{1-a}
\frac{(a;q)_{2k}}{(q;q)_{2k}}\frac{(b,q/b;q)_k}{(aq^{3n+1},q^{-3n};q)_k}\frac{(aq^{3n+1},q^{-3n};q^3)_k}{(aq^{3}/b,abq^{2};q^3)_k}q^kE_k(a,b)
\\[3pt]
&\quad=\frac{(aq;q)_{3n}}{(q;q)_{3n}}\frac{(q^3,bq,q^{2}/b;q^3)_n}{(aq,aq^{3}/b,abq^{2};q^3)_n}F_n(a,b),
\end{align*}
where
\begin{align*}
&E_k(a,b)=\sum_{i=1}^k\frac{q^{i-1}}{1-bq^{i-1}}-\sum_{i=1}^k\frac{q^{i}/b^2}{1-q^{i}/b}
+\sum_{i=1}^k\frac{aq^{3i}/b^2}{1-aq^{3i}/b}-\sum_{i=1}^k\frac{aq^{3i-1}}{1-abq^{3i-1}},
\\[2pt]
&F_n(a,b)=\sum_{j=1}^n\frac{q^{3j-2}}{1-bq^{3j-2}}-\sum_{j=1}^n\frac{q^{3j-1}/b^2}{1-q^{3j-1}/b}
+\sum_{j=1}^n\frac{aq^{3j}/b^2}{1-aq^{3j}/b}-\sum_{j=1}^n\frac{aq^{3j-1}}{1-abq^{3j-1}}.
\end{align*}
Through the operator $\mathcal{D}_{b}$ and the last equation, it is
clear that
\begin{align}
&\sum_{k=1}^{n}\frac{1-aq^{4k}}{1-a}
\frac{(a;q)_{2k}}{(q;q)_{2k}}\frac{(b,q/b;q)_k}{(aq^{3n+1},q^{-3n};q)_k}\frac{(aq^{3n+1},q^{-3n};q^3)_k}{(aq^{3}/b,abq^{2};q^3)_k}q^k
\bigg\{E_k(a,b)^2-G_k(a,b)\bigg\}
\notag\\[3pt]
&\quad=\frac{(aq;q)_{3n}}{(q;q)_{3n}}\frac{(q^3,bq,q^{2}/b;q^3)_n}{(aq,aq^{3}/b,abq^{2};q^3)_n}\bigg\{F_n(a,b)^2-H_n(a,b)\bigg\},
\label{eq:wei-bbb}
\end{align}
where
\begin{align*}
&G_k(a,b)=\sum_{i=1}^k\frac{q^{2i-2}}{(1-bq^{i-1})^2}-\sum_{i=1}^k\frac{(q^{i}/b-2)q^{i}/b^3}{(1-q^{i}/b)^2}
\end{align*}
\begin{align*}
&\qquad\qquad\:+\sum_{i=1}^k\frac{(aq^{3i}/b-2)aq^{3i}/b^3}{(1-aq^{3i}/b)^2}-\sum_{i=1}^k\frac{a^2q^{6i-2}}{(1-abq^{3i-1})^2},
\\[2pt]
&H_n(a,b)=\sum_{j=1}^n\frac{q^{6j-4}}{(1-bq^{3j-2})^2}-\sum_{j=1}^n\frac{(q^{3j-1}/b-2)q^{3j-1}/b^3}{(1-q^{3j-1}/b)^2}
\\[2pt]
&\qquad\qquad\:+\sum_{j=1}^n\frac{(aq^{3j}/b-2)aq^{3j}/b^3}{(1-aq^{3j}/b)^2}-\sum_{j=1}^n\frac{a^2q^{6j-2}}{(1-abq^{3j-1})^2}.
\end{align*}
The $a\to q, b\to q, q\to q^2$ case of \eqref{eq:wei-bbb} can be
manipulated as
\begin{align*}
&\sum_{k=1}^{n}q^{2k}[8k+1]\frac{(q;q^2)_{k}^2(q;q^2)_{2k}}{(q^2;q^2)_{2k}(q^6;q^6)_{k}^2}\frac{(q^{3+6n},q^{-6n};q^6)_k}{(q^{3+6n},q^{-6n};q^2)_k}
\sum_{i=1}^{k}\bigg\{\frac{q^{2i-1}}{[2i-1]^2}-\frac{q^{6i}}{[6i]^2}\bigg\}
\notag\\[3pt]
&\quad=\frac{(q^3;q^2)_{3n}(q^3;q^6)_{n}}{(q^2;q^2)_{3n}(q^6;q^6)_{n}}
\sum_{j=1}^{2n}(-1)^{j-1}\frac{q^{3j}}{[3j]^2}.
\end{align*}
Letting $n\to \infty$ in the upper identity, we are led to Theorem
\ref{thm-e}.
\end{proof}

For the aim to prove Theorem \ref{thm-f}, we shall draw support from
the summation formula for basic hypergeometric series  (cf. \cite[p.
65]{Chu-a}):
\begin{align}
&\sum_{k=0}^{n}\frac{1-a^{\frac{2}{3}}q^{\frac{4k}{3}}}{1-a^{\frac{2}{3}}}
\frac{(q^{-n},a^{\frac{2}{3}}q^{n+\frac{2}{3}};q)_k}{(qa^{\frac{2}{3}}b^{-\frac{1}{3}},q^{\frac{1}{3}}a^{\frac{2}{3}}b^{\frac{1}{3}};q)_k}
\frac{(b^{\frac{1}{3}},q^{\frac{2}{3}}b^{-\frac{1}{3}};q^{\frac{1}{3}})_{k}}{(a^{\frac{2}{3}}q^{n+\frac{1}{3}},q^{-n-\frac{1}{3}};q^{\frac{1}{3}})_{k}}
\frac{(a^{\frac{2}{3}};q^{\frac{1}{3}})_{2k}}{(q^{\frac{1}{3}};q^{\frac{1}{3}})_{2k}}q^{\frac{k}{3}}
\notag\\[3pt]
  &\quad
=\frac{(q^{\frac{2}{3}}a^{\frac{2}{3}};q^{\frac{1}{3}})_{n}}{(q^{\frac{1}{3}};q^{\frac{1}{3}})_{n}}
\frac{(q,q^{\frac{2}{3}}b^{\frac{1}{3}},q^{\frac{4}{3}}b^{-\frac{1}{3}};q)_n}
{(q^{\frac{2}{3}}a^{\frac{2}{3}},qa^{\frac{2}{3}}b^{-\frac{1}{3}},q^{\frac{1}{3}}a^{\frac{2}{3}}b^{\frac{1}{3}};q)_n},
\label{basic-c}
\end{align}
 where we have replaced
$$\frac{(b^{\frac{1}{3}},q^{\frac{2}{3}}b^{-\frac{1}{3}};q)_{k}}{(a^{\frac{2}{3}}q^{n+\frac{1}{3}},q^{-n-\frac{1}{3}};q)_{k}}
\quad\text{by}\quad
\frac{(b^{\frac{1}{3}},q^{\frac{2}{3}}b^{-\frac{1}{3}};q^{\frac{1}{3}})_{k}}{(a^{\frac{2}{3}}q^{n+\frac{1}{3}},q^{-n-\frac{1}{3}};q^{\frac{1}{3}})_{k}}$$
 for correction. The nonterminating form of
\eqref{basic-c} can be observed in Gasper and Rahman \cite[Equation
(4.5)]{Gasper-Rahman}.

Finally, we begin to prove Theorem \ref{thm-f}.

\begin{proof}[{\bf{Proof of Theorem \ref{thm-f}}}]
Apply the partial derivative operator $\mathcal{D}_{b}$ to the $a\to
a^{3/2}, b\to b^3, q\to q^3$ case of \eqref{basic-c} to deduce
\begin{align*}
&\sum_{k=1}^{n}\frac{1-aq^{4k}}{1-a}
\frac{(a/q;q)_{2k}}{(q^2;q)_{2k}}\frac{(b,q^2/b;q)_k}{(aq^{3n+1},q^{-1-3n};q)_k}\frac{(aq^{3n+2},q^{-3n};q^3)_k}{(aq^{3}/b,abq;q^3)_k}q^kR_k(a,b)
\\[3pt]
&\quad=\frac{(aq;q)_{3n}}{(q^2;q)_{3n}}\frac{(q^3,bq^2,q^{4}/b;q^3)_n}{(aq^2,aq^{3}/b,abq;q^3)_n}S_n(a,b),
\end{align*}
where
\begin{align*}
&R_k(a,b)=\sum_{i=1}^k\frac{q^{i-1}}{1-bq^{i-1}}-\sum_{i=1}^k\frac{q^{i+1}/b^2}{1-q^{i+1}/b}
+\sum_{i=1}^k\frac{aq^{3i}/b^2}{1-aq^{3i}/b}-\sum_{i=1}^k\frac{aq^{3i-2}}{1-abq^{3i-2}},
\\[2pt]
&S_n(a,b)=\sum_{j=1}^n\frac{q^{3j-1}}{1-bq^{3j-1}}-\sum_{j=1}^n\frac{q^{3j+1}/b^2}{1-q^{3j+1}/b}
+\sum_{j=1}^n\frac{aq^{3j}/b^2}{1-aq^{3j}/b}-\sum_{j=1}^n\frac{aq^{3j-2}}{1-abq^{3j-2}}.
\end{align*}
Employing the operator $\mathcal{D}_{b}$ to both sides of the last
equation, it is obvious that
\begin{align}
&\sum_{k=1}^{n}\frac{1-aq^{4k}}{1-a}
\frac{(a/q;q)_{2k}}{(q^2;q)_{2k}}\frac{(b,q^2/b;q)_k}{(aq^{3n+1},q^{-1-3n};q)_k}\frac{(aq^{3n+2},q^{-3n};q^3)_k}{(aq^{3}/b,abq;q^3)_k}q^k
\bigg\{R_k(a,b)^2-U_k(a,b)\bigg\}
\notag\\[3pt]
&\quad=\frac{(aq;q)_{3n}}{(q^2;q)_{3n}}\frac{(q^3,bq^2,q^{4}/b;q^3)_n}{(aq^2,aq^{3}/b,abq;q^3)_n}\bigg\{S_n(a,b)^2-V_n(a,b)\bigg\},
\label{eq:wei-ccc}
\end{align}
where
\begin{align*}
&U_k(a,b)=\sum_{i=1}^k\frac{q^{2i-2}}{(1-bq^{i-1})^2}-\sum_{i=1}^k\frac{(q^{i+1}/b-2)q^{i+1}/b^3}{(1-q^{i+1}/b)^2}
\\[2pt]
&\qquad\qquad\:+\sum_{i=1}^k\frac{(aq^{3i}/b-2)aq^{3i}/b^3}{(1-aq^{3i}/b)^2}-\sum_{i=1}^k\frac{a^2q^{6i-4}}{(1-abq^{3i-2})^2},
\\[2pt]
&V_n(a,b)=\sum_{j=1}^n\frac{q^{6j-2}}{(1-bq^{3j-1})^2}-\sum_{j=1}^n\frac{(q^{3j+1}/b-2)q^{3j+1}/b^3}{(1-q^{3j+1}/b)^2}
\\[2pt]
&\qquad\qquad\:+\sum_{j=1}^n\frac{(aq^{3j}/b-2)aq^{3j}/b^3}{(1-aq^{3j}/b)^2}-\sum_{j=1}^n\frac{a^2q^{6j-4}}{(1-abq^{3j-2})^2}.
\end{align*}
The $a\to q^{-1}, b\to q^2, q\to q^2$ case of \eqref{eq:wei-ccc}
engenders
\begin{align*}
&\sum_{k=1}^{n}q^{2k}[8k-1]\frac{(q^2;q^2)_{k}^2(q^{-3};q^2)_{2k}}{(q^4;q^2)_{2k}(q^3;q^6)_{k}^2}\frac{(q^{3+6n},q^{-6n};q^6)_k}{(q^{1+6n},q^{-2-6n};q^2)_k}
\sum_{i=1}^{k}\bigg\{\frac{q^{6i-3}}{[6i-3]^2}-\frac{q^{2i}}{[2i]^2}\bigg\}
\notag\\[3pt]
&\quad=\frac{(q;q^2)_{3n}(q^6;q^6)_{n}^3}{(q^4;q^2)_{3n}(q^3;q^6)_{n}^3}
\sum_{j=1}^{2n}(-1)^{j}\frac{q^{3j-1}}{[3j]^2}.
\end{align*}
Letting $n\to \infty$ in this identity, we catch hold of Theorem
\ref{thm-f}.
\end{proof}

\textbf{Acknowledgments}

 The work is supported by the National Natural Science Foundation of China (No. 12071103).


\end{document}